\documentclass[12pt]{article}
\usepackage{times}
\usepackage{amsmath}
\usepackage{amssymb}
\usepackage{amsthm}
\usepackage{amsmath}
\usepackage{amsfonts}
\usepackage{amssymb}
\usepackage[all]{xy}
\usepackage{xcolor}
\usepackage[utf8]{inputenc}

\newtheorem{theorem}{Theorem}[section]

\newtheorem{corollary}[theorem]{Corollary}
\newtheorem{proposition}[theorem]{Proposition}

\theoremstyle{definition}

\newtheorem{definition}[theorem]{Definition}

\def\0{\underline 0}

\begin{document}

\title {\bf About relative polar varieties and Brasselet numbers \footnote{
\newline   {\it } }}

\vspace{1cm}
\author{Hellen Santana\\Universidade de São Paulo}
\date{}

\maketitle

\begin{abstract}
\noindent In this work, we study the consequences of an empty polar variety on the topology of a function-germ with (possibly) nonisolated singularities defined on a singular variety.

\end{abstract}

\section*{Introduction}
\hspace{0,5cm} Let $f:(\mathbb{C}^n,0)\rightarrow(\mathbb{C},0)$ be an analytic function defined in a neighbourhood of the origin. In \cite{Milnor}, Milnor described the topology of the germ $f$ at the origin, proving that,  if $f$ has an isolated singularity at the origin, the set $f^{-1}(\delta)\cap B_{\epsilon}$, later called Milnor fiber, has the homotopy type of a bouquet of $\mu(f)$ spheres of dimension $n-1,$ where $\delta$ is a regular value of $f, 0<|\delta|\ll\epsilon\ll1$. The number $\mu(f)$ is the Milnor number of $f$, an invariant associated to the germ $f$ which also counts the number of Morse points in a Morsefication of $f$ in a neighbourhood of the origin. 

The Milnor number was largely studied and generalized to many settings (for instance, for mention a few \cite{T}, \cite{Hamm}, \cite{L},\cite{Greuel}, \cite{massey2003numerical},\cite{iomdin1974complex}). If $f$ is defined over a complex analytic space $X$ and $f$ has an isolated singularity at the origin, a generalization for the Milnor number is the Euler obstruction of the function $f$, introduced in \cite{BMPS}, by Brasselet, Massey, Parameswaran and Seade. In \cite{STV}, Seade, Tib\u{a}r and Verjovsky proved that, up to sign, this number is the number of Morse critical points of a stratified Morsefication of $f$ appearing in the regular part of $X$ in a neighbourhood of the origin.   


In the case where $f$ is defined over a complex analytic germ $(X,0)$ equipped with a good stratification $\mathcal{V}$ relative to $f$ and the function $f$ has nonisolated singularities, Dutertre and Grulha provided a way to numerically describe the topology of the generalized Milnor fiber $X\cap f^{-1}(\delta)\cap B_{\epsilon}$. In \cite{DG}, the authors defined the Brasselet number $B_{f,X}(0)$ of $f$ at the origin and proved a Lê-Greuel type formula for this number: if $g:X\rightarrow\mathbb{C}$ is prepolar with respect to $\mathcal{V}$ at the origin and $0<|\delta|\ll\epsilon\ll1,$ then $B_{f,X}(0)-B_{f,X^g}(0)=(-1)^{d-1}n_q,$ where $n_q$ is the number of Morse critical points of a partial Morsefication of $g|_{X\cap f^{-1}(\delta)\cap B_{\epsilon}}$ appearing in the regular part of $X$, and $X^g=X\cap \{g=0\}.$

Computing these number of stratified Morse critical points is directly connected to relative polar varieties. Consider a linear form $l$ in $\mathbb{C}^n$, a Whitney stratification of $X$ and a function-germ $g:(U,0)\rightarrow(\mathbb{C},0)$. If $l$ is sufficiently generic, the polar variety (curve) $\Gamma_{g,l}$ defined by Lê and Teissier, in \cite{LT}, coincides with the relative polar curve defined by Massey, in \cite{Ms1}, and with the relative polar varieties defined by Massey in \cite{massey2003numerical}(see \cite{massey2008enriched} and \cite{massey2018IPA}). Each of these polar varieties is useful to compute polar multiplicities (\cite{loeser1984formules},\cite{teissier1982varietes}) and intersection numbers (\cite{massey2003numerical}), but also to describe the critical loci of a pair of functions defined over $X$ (\cite{DG},\cite{Ms1}), which is the approach we are interested the most. 
In this work, we use polar varieties to compute a number of stratified Morse critical points of a specific type of deformation of a function-germ aiming to obtain informations about the Brasselet number of this germ.

\section{Local Euler obstruction and Euler obstruction of a function}

\hspace{0,5cm} In this section, we will see the definition of the local Euler obstruction, a singular invariant defined by MacPherson and used as one of the main tools in his proof of the Deligne-Grothendieck conjecture about the existence and uniqueness of Chern classes for singular varities (see \cite{MacPherson}). 

Let $(X,0)\subset(\mathbb{C}^n,0)$ be an equidimensional reduced complex analytic germ of dimension $d$ in a open set $U\subset\mathbb{C}^n.$ Consider a complex analytic Whitney stratification $\mathcal{V}=\{V_{\lambda}\}$ of $U$ adapted to $X$ such that $\{0\}$ is a stratum. We choose a small representative of $(X,0),$ denoted by $X,$ such that $0$ belongs to the closure of all strata. We write $X=\cup_{i=0}^{q} V_i,$ where $V_0=\{0\}$ and $V_q=X_{reg},$ where $X_{reg}$ is the regular part of $X.$ We suppose that $V_0,V_1,\ldots,V_{q-1}$ are connected and that the analytic sets $\overline{V_0},\overline{V_1},\ldots,\overline{V_q}$ are reduced. We write $d_i=dim(V_i), \ i\in\{1,\ldots,q\}.$ Note that $d_q=d.$ 

Let $G(d,N)$ be the Grassmannian manifold, $x\in X_{reg}$ and consider the Gauss map $\phi: X_{reg}\rightarrow U\times G(d,N)$ given by $x\mapsto(x,T_x(X_{reg})).$ 

\begin{definition}
The closure of the image of the Gauss map $\phi$ in $U\times G(d,N)$, denoted by $\tilde{X}$, is called \textbf{Nash modification} of $X$. It is a complex analytic space endowed with an analytic projection map $\nu:\tilde{X}\rightarrow X.$
\end{definition}

Consider the extension of the tautological bundle $\mathcal{T}$ over $U\times G(d,N).$ Since \linebreak$\tilde{X}\subset U\times G(d,N)$, we consider $\tilde{T}$ the restriction of $\mathcal{T}$ to $\tilde{X},$ called the \textbf{Nash bundle}, and $\pi:\tilde{T}\rightarrow\tilde{X}$ the projection of this bundle.

In this context, denoting by $\varphi$ the natural projection of $U\times G(d,N)$ at $U,$ we have the following diagram:

$$\xymatrix{
\tilde{T} \ar[d]_{\pi}\ar[r] & \mathcal{T}\ar[d] \\ 
\tilde{X}\ar[d]_{\nu}\ar[r] & U\times G(d,N)\ar[d]^{\varphi} \\ 
X\ar[r] & U\subseteq\mathbb{C}^N \\}  $$

Considering $\vert\vert z\vert\vert=\sqrt{z_1\overline{z_1}+\cdots+z_N\overline{z_N}}$, the $1$-differential form $w=d\vert\vert z\vert\vert^2$ over $\mathbb{C}^N$ defines a section in $T^{*}\mathbb{C}^N$ and its pullback $\varphi^{*}w$ is a $1$- form over $U\times G(d,N).$ Denote by $\tilde{w}$ the restriction of $\varphi^{*}w$ over $\tilde{X}$, which is a section of the dual bundle $\tilde{T}^{*}.$

Choose $\epsilon$ small enough for $\tilde{w}$ be a non zero section over $\nu^{-1}(z), 0<\vert\vert z \vert\vert\leqslant\epsilon,$ let $B_{\epsilon}$ be the closed ball with center at the origin with radius $\epsilon$ and denote by:

\begin{enumerate}

\item $Obs(\tilde{T}^{*},\tilde{w})\in\mathbb{H}^{2d}(\nu^{-1}(B_{\epsilon}),\nu^{-1}(S_{\epsilon}),\mathbb{Z})$  the obstruction for extending $\tilde{w}$ from $\nu^{-1}(S_{\epsilon})$ to $\nu^{-1}(B_{\epsilon});$

\item $O_{\nu^{-1}(B_{\epsilon}),\nu^{-1}(S_{\epsilon})}$ the fundamental class in $\mathbb{H}_{2d}(\nu^{-1}(B_{\epsilon}),\nu^{-1}(S_{\epsilon}),\mathbb{Z}).$ 
\end{enumerate}

\begin{definition}
The \textbf{local Euler obstruction} of $X$ at $0, \ Eu_X(0),$ is given by the evaluation $$Eu_X(0)=\langle Obs(\tilde{T}^{*},\tilde{w}),O_{\nu^{-1}(B_{\epsilon}),\nu^{-1}(S_{\epsilon})}\rangle.$$
\end{definition}

In \cite{BLS}, Brasselet, Lê and Seade proved a formula to make the calculation of the Euler obstruction easier.

\begin{theorem}(Theorem 3.1 of \cite{BLS})
Let $(X,0)$ and $\mathcal{V}$ be given as before, then for each generic linear form $l,$ there exists $\epsilon_0$ such that for any $\epsilon$ with $0<\epsilon<\epsilon_0$ and $\delta\neq0$ sufficiently small, the Euler obstruction of $(X,0)$ is equal to 

$$Eu_X(0)=\sum^{q}_{i=1}\chi(V_i\cap B_{\epsilon}\cap l^{-1}(\delta)).Eu_{X}(V_i),$$

\noindent where $\chi$ is the Euler characteristic, $Eu_{X}(V_i)$ is the Euler obstruction of $X$ at a point of $V_i, \ i=1,\ldots,q$ and $0<|\delta|\ll\epsilon\ll1.$
\end{theorem} 


Let us give the definition of another invariant introduced by Brasselet, Massey, Parameswaran and Seade in \cite{BMPS}. Let $f:X\rightarrow\mathbb{C}$ be a holomorphic function with isolated singularity at the origin given by the restriction of a holomorphic function $F:U\rightarrow\mathbb{C}$ and denote by $\overline{\nabla}F(x)$ the conjugate of the gradient vector field of $F$ in $x\in U,$ $$\overline{\nabla}F(x):=\left(\overline{\frac{\partial F}{\partial x_1}},\ldots, \overline{\frac{\partial F}{\partial x_n}}\right).$$

Since $f$ has an isolated singularity at the origin, for all $x\in X\setminus\{0\},$ the projection $\hat{\zeta}_i(x)$ of $\overline{\nabla}F(x)$ over $T_x(V_i(x))$ is nonzero, where $V_i(x)$ is a stratum containing $x.$ Using this projection, the authors constructed, in \cite{BMPS}, a stratified vector field over $X,$ denoted by $\overline{\nabla}f(x).$ Let $\tilde{\zeta}$ be the lifting of $\overline{\nabla}f(x)$ as a section of the Nash bundle $\tilde{T}$ over $\tilde{X}$, without singularity over $\nu^{-1}(X\cap S_{\epsilon}).$

Let $\mathcal{O}(\tilde{\zeta})\in\mathbb{H}^{2n}(\nu^{-1}(X\cap B_{\epsilon}),\nu^{-1}(X\cap S_{\epsilon}))$ be the obstruction cocycle for extending $\tilde{\zeta}$ as a non zero section of $\tilde{T}$ inside $\nu^{-1}(X\cap B_{\epsilon}).$

\begin{definition}
The \textbf{local Euler obstruction of the function} $f, Eu_{f,X}(0)$ is the evaluation of $\mathcal{O}(\tilde{\zeta})$ on the fundamental class $[\nu^{-1}(X\cap B_{\epsilon}),\nu^{-1}(X\cap S_{\epsilon})].$
\end{definition}

The next theorem compares the Euler obstruction of a space $X$ with the Euler obstruction of function defined over $X.$

\begin{theorem}\label{Euler obstruction of a function formula}(Theorem 3.1 of \cite{BMPS})
Let $(X,0)$ and $\mathcal{V}$ be given as before and let \linebreak$f:(X,0)\rightarrow(\mathbb{C},0)$ be a function with an isolated singularity at $0.$ For $0<|\delta|\ll\epsilon\ll1,$ we have
 $$Eu_{f,X}(0)=Eu_X(0)-\sum_{i=1}^{q}\chi(V_i\cap B_{\epsilon}\cap f^{-1}(\delta)).Eu_X(V_i).$$
\end{theorem}



Let us now see a definition we will need to define a generic point of a function-germ. Let $\mathcal{V}=\{V_{\lambda}\}$ be a stratification of a reduced complex analytic space $X.$

\begin{definition}
Let $p$ be a point in a stratum $V_{\beta}$ of $\mathcal{V}.$ A \textbf{degenerate tangent plane of $\mathcal{V}$ at $p$} is an element $T$ of some Grassmanian manifold such that $T=\displaystyle\lim_{p_i\rightarrow p}T_{p_i}V_{\alpha},$ where $p_i\in V_{\alpha}$, $V_{\alpha}\neq V_{\beta}.$
\end{definition}

\begin{definition}
Let $(X,0)\subset(U,0)$ be a germ of complex analytic space in $\mathbb{C}^n$ equipped with a Whitney stratification and let $f:(X,0)\rightarrow(\mathbb{C},0)$ be an analytic function, given by the restriction of an analytic function $F:(U,0)\rightarrow(\mathbb{C},0).$ Then $0$ is said to be a \textbf{generic point}\index{holomorphic function germ!generic point of} of $f$ if the hyperplane $Ker(d_0F)$ is transverse in $\mathbb{C}^n$ to all degenerate tangent planes of the Whitney stratification at $0.$ 
\end{definition}

Now, let us see the definition of a Morsefication of a function. 

\begin{definition}
Let $\mathcal{W}=\{W_0,W_1,\ldots,W_q\},$ with $0\in W_0,$ a Whitney stratification of the complex analytic space $X.$ A function $f:(X,0)\rightarrow(\mathbb{C},0)$ is said to be \textbf{Morse stratified} if $\dim W_0\geq1, f|_{W_0}: W_0\rightarrow\mathbb{C}$ has a Morse point at $0$ and $0$ is a generic point of $f$ with respect to $W_{i},$ for all $ i\neq0.$
\end{definition}

A \textbf{stratified Morsefication}\index{holomorphic function germ!stratified Morsefication of} of a germ of analytic function $f:(X,0)\rightarrow(\mathbb{C},0)$ is a deformation $\tilde{f}$ of $f$ such that $\tilde{f}$ is Morse stratified.

In \cite{STV}, Seade, Tib\u{a}r and Verjovsky proved that the Euler obstruction of a function $f$ is also related to the number of Morse critical points of a stratified Morsefication of $f.$

\begin{proposition}(Proposition 2.3 of \cite{STV})\label{Eu_f and Morse points}
Let $f:(X,0)\rightarrow(\mathbb{C},0)$ be a germ of analytic function with isolated singularity at the origin. Then, \begin{center}
$Eu_{f,X}(0)=(-1)^dn_{reg},$
\end{center}
where $n_{reg}$ is the number of Morse points in $X_{reg}$ in a stratified Morsefication of $f.$
\end{proposition}

\section{Brasselet number}

\hspace{0,5cm} In this section, we present definitions and results needed in the development of the results of this work. The main reference for this section is \cite{Ms1}.

Let $X$ be a reduced complex analytic space (not necessarily equidimensional) of dimension $d$ in an open set $U\subseteq\mathbb{C}^n$ and let $f:(X,0)\rightarrow(\mathbb{C},0)$ be an analytic map. We write $V(f)=f^{-1}(0).$ 

\begin{definition}\label{good stratification}
A \textbf{good stratification of $X$ relative to $f$} is a stratification $\mathcal{V}$ of $X$ which is adapted to $V(f)$ such that $\{V_{\lambda}\in\mathcal{V},V_{\lambda}\nsubseteq V(f)\}$ is a Whitney stratification of $X\setminus V(f)$ and such that for any pair $(V_{\lambda},V_{\gamma})$ such that $V_{\lambda}\nsubseteq V(f)$ and $V_{\gamma}\subseteq V(f),$ the $(a_f)$-Thom condition is satisfied, that is, if $p\in V_{\gamma}$ and $p_i\in V_{\lambda}$ are such that $p_i\rightarrow p$ and $T_{p_i} V(f|_{V_{\lambda}}-f|_{V_{\lambda}}(p_i))$ converges to some $\mathcal{T},$ then $T_p V_{\gamma}\subseteq\mathcal{T}.$
\end{definition}

If $f:X\rightarrow\mathbb{C}$ has a stratified isolated critical point and $\mathcal{V}$ is a Whitney stratification of $X,$ then \begin{equation}\label{induced stratification}
    \{V_{\lambda}\setminus X^f, V_{\lambda}\cap X^f\setminus\{0\},\{0\}, V_{\lambda}\in\mathcal{V}\}
\end{equation}

\noindent is a good stratification of $X$ relative to $f,$ called the good stratification induced by $f.$

Let $\mathcal{V}$ be a good stratification of $X$ relative to $f.$

\begin{definition}
The \textbf{critical locus of $f$ relative to $\mathcal{V}$}, $\Sigma_{\mathcal{V}}f,$ is given by the union \begin{center}$\Sigma_{\mathcal{V}}f=\displaystyle\bigcup_{V_{\lambda}\in\mathcal{V}}\Sigma(f|_{V_{\lambda}}).$\end{center}
\end{definition}

\begin{definition}
If $\mathcal{V}=\{V_{\lambda}\}$ is a stratification of $X,$ the \textbf{relative polar variety of $f$ and $g$ with respect to $\mathcal{V}$}, denoted by $\Gamma_{f,g}(\mathcal{V}),$ is the the union $\cup_{\lambda}\Gamma_{f,g}(V_{\lambda}),$ where $\Gamma_{f,g}(V_{\lambda})$ denotes the closure in $X$ of the critical locus of $(f,g)|_{V_{\lambda}\setminus X^f},$ where $X^f=X\cap\{f=0\}.$  
\end{definition}

\begin{definition}
If $\mathcal{V}=\{V_{\lambda}\}$ is a stratification of $X,$ the \textbf{symmetric relative polar variety of $f$ and $g$ with respect to $\mathcal{V}$}, $\tilde{\Gamma}_{f,g}(\mathcal{V}),$ is the union $\cup_{\lambda}\tilde{\Gamma}_{f,g}(V_{\lambda}),$ where $\Gamma_{f,g}(V_{\lambda})$ denotes the closure in $X$ of the critical locus of $(f,g)|_{V_{\lambda}\setminus (X^f\cup X^g)},$  $X^f=X\cap \{f=0\}$ and $X^g=X\cap \{g=0\}. $ 
\end{definition}

\begin{definition}\label{definition prepolar}
Let $\mathcal{V}$ be a good stratification of $X$ relative to a function
$f:(X,0)\rightarrow(\mathbb{C},0).$ A function $g :(X, 0)\rightarrow(\mathbb{C},0)$ is \textbf{prepolar with respect to $\mathcal{V}$ at the origin} if the origin is a stratified isolated critical point, that is, $0$ is an isolated point of $\Sigma_{\mathcal{V}}g.$
\end{definition} 




\begin{definition}\label{definition tractable}
A function $g :(X, 0)\rightarrow(\mathbb{C},0)$ is \textbf{tractable at the origin with respect to a good stratification $\mathcal{V}$ of $X$ relative to $f :(X, 0)\rightarrow(\mathbb{C},0)$} if $dim_0 \ \tilde{\Gamma}^1_{f,g}(\mathcal{V})\leq1$ and, for all strata $V_{\alpha}\subseteq X^f$,
$g|_{V_{\alpha}}$ has no critical point in a neighbourhood of the origin except perhaps at the origin itself.

\end{definition}

Another concept useful for this work is the notion of constructible functions. Consider a Whitney stratification $\mathcal{W}=\{W_1,\ldots, W_q\}$ of $X$ such that each stratum $W_i$ is connected.

\begin{definition}
A constructible function with respect to the stratification $\mathcal{W}$ of $X$ is a function $\beta:X\rightarrow\mathbb{Z}$ which is constant on each stratum $W_i,$ that is, there exist integers $t_1,\ldots,t_q,$ such that $\beta=\sum_{i=1}^{q}t_i.1_{W_i},$ where $1_{W_i}$ is the characteristic function of $W_i.$
\end{definition}

\begin{definition}
The Euler characteristic $\chi(X,\beta)$ of a constructible function $\beta:X\rightarrow\mathbb{Z}$ with respect to the stratification $\mathcal{W}$ of $X,$ given by $\beta=\sum_{i=1}^{q}t_i.1_{W_i}$, is defined by $\chi(X,\beta)=\sum_{i=1}^{q}t_i.\chi(W_i).$
\end{definition}

We present now the definition of the Brasselet number and the main theorems of \cite{DG}, used as inspiration for this work.

Let $f: (X,0)\rightarrow(\mathbb{C},0)$ be a complex analytic function germ and let $\mathcal{V}$ be a good stratification of $X$ relative to $f.$ We denote by $V_1,\ldots, V_q$ the strata of $\mathcal{V}$ that are not contained in $\{f=0\}$ and we assume that $V_1,\ldots, V_{q-1}$ are connected and that $V_{q}=\linebreak X_{reg}\setminus \{f=0\}.$ Note that $V_q$ could be not connected.  
 
\begin{definition}
Suppose that $X$ is equidimensional. Let $\mathcal{V}$ be a good stratification of $X$ relative to $f.$ The \textbf{Brasselet number} of $f$ at the origin, $B_{f,X}(0),$ is defined by \begin{center}
$B_{f,X}(0)=\sum_{i=1}^{q}\chi(V_i\cap f^{-1}(\delta)\cap B_{\epsilon})Eu_X(V_i),$
\end{center}
where $0<|\delta|\ll\epsilon\ll1.$
\end{definition} 

\noindent\textbf{Remark:} If $V_q^i$ is a connected component of $V_{q},$ $Eu_X(V_q^i)=1.$

Notice that if $f$ has a stratified isolated singularity at the origin, then \linebreak$B_{f,X}(0)=Eu_{X}(0)-Eu_{f,X}(0)$ (see Theorem \ref{Euler obstruction of a function formula}).

In \cite{DG}, Dutertre and Grulha proved interesting formulas describing the topological relation between the Brasselet number and a number of certain critical points of a special type of deformation of functions. Let us now present some of these results. Fist we need the definition of a special type of Morsefication, introduced by Dutertre and Grulha.

\begin{definition}
A \textbf{partial Morsefication} of $g:f^{-1}(\delta)\cap X\cap B_{\epsilon}\rightarrow\mathbb{C}$ is a function $\tilde{g}: f^{-1}(\delta)\cap X\cap B_{\epsilon}\rightarrow\mathbb{C}$ (not necessarily holomorphic) which is a local Morsefication of all isolated critical points of $g$ in $f^{-1}(\delta)\cap X\cap \{g\neq 0\}\cap B_{\epsilon}$ and which coincides with $g$ outside a small neighbourhood of these critical points.
\end{definition} 

Let $g : (X, 0) \rightarrow (\mathbb{C},0)$ be a complex analytic function which is tractable at the origin with respect to $\mathcal{V}$ relative to $f .$ Then $\tilde{\Gamma}_{f,g}$ is a complex analytic curve and for $0<|\delta|\ll1$ the critical points of $g|_{f^{-1}(\delta)\cap X}$ in $B_{\epsilon}$ lying outside $\{g=0\}$ are isolated.
Let $\tilde{g}$ be a partial Morsefication of $g:f^{-1}(\delta)\cap X\cap B_{\epsilon}\rightarrow\mathbb{C}$ and, for each $i\in\{1,\ldots, q\},$ let $n_q$ be the number of stratified Morse critical points of $\tilde{g}$ appearing on $X_{reg}\cap f^{-1}(\delta)\cap \{g\neq 0\}\cap B_{\epsilon}.$ 

\begin{theorem}\label{4.4 DG}(Theorem 4.4 of \cite{DG})
Suppose that $X$ is equidimensional and that $g$ is prepolar with respect to $\mathcal{V}$ at the origin. For $0<|\delta|\ll\epsilon\ll 1,$ we have\begin{center}
 $B_{f,X}(0)-B_{f,X^g}(0)=(-1)^{d-1}n_q,$
 \end{center}
 where $n_q$ is the number of stratified Morse critical points on the top stratum $V_q\cap f^{-1}(\delta)\cap B_{\epsilon}$ appearing in a Morsefication of $g:X\cap f^{-1}(\delta)\cap B_{\epsilon}\rightarrow \mathbb{C}.$
\end{theorem}


In the case where $X$ is equipped with a Whitney stratification $\mathcal{V}=\{V_0,V_1,\ldots,V_q\}$ with $V_0=\{0\},$ and $f,g:X\rightarrow\mathbb{C}$ have an isolated stratified singularity at the origin with respect to this stratification. In \cite{DG}, the authors also related the topology of the  generalized Minor fibres of $f$ and $g$ and some number of Morse points. 

\begin{corollary}\label{corollary 6.5 of DG}
Suppose that $X$ is equidimensional and that $g$ (resp. $f$) is prepolar with respect to the good stratification induced by $f$ (resp. $g$) at the origin. Then \begin{center}
$B_{f,X}(0)-B_{g,X}(0)=(-1)^{d-1}(n_q-m_q),$
\end{center}
\noindent where $n_q$ (resp. $m_q$) is the number of stratified Morse critical points on the top stratum \linebreak$V_q\cap f^{-1}(\delta)\cap B_{\epsilon}$ (resp. $V_q\cap g^{-1}(\delta)\cap B_{\epsilon}$) appearing in a Morsefication of \linebreak$g:X\cap f^{-1}(\delta)\cap B_{\epsilon}\rightarrow\mathbb{C}$ (resp. $f:X\cap g^{-1}(\delta)\cap B_{\epsilon}\rightarrow\mathbb{C}$).
\end{corollary}

\section{Brasselet numbers and empty polar varieties}

In this final section, we present relations between Brasselet numbers of two function-germs in the case where the relative polar variety associated to these germs is empty.
Let $(X,0)$ be a complex analytic space in $U\subset\mathbb{C}^{n+1}$. Let $f,g:(U,0)\rightarrow(\mathbb{C},0)$ be germs of holomorphic functions and let $\mathcal{V}$ be a good stratification of $X$ relative to $f.$ Suppose that $\Sigma_{\mathcal{V}}g\cap\{f=0\}=\{0\}$. We aim to obtain information about the Brasselet number of $f$ and the Brasselet number of $g$ in the case where the relative polar variety $\Gamma_{f,g}(\mathcal{V})$ is empty. 
We begin with a description of two relevant subsets of $\Gamma_{f,g}(\mathcal{V})$.

\begin{proposition}\label{3.1}
The stratified critical set $\Sigma_{\mathcal{V}}g$ of $g$ and the symmetric relative polar variety $\tilde{\Gamma}_{f,g}(\mathcal{V})$ are subsets of $\Gamma_{f,g}(\mathcal{V}).$
\end{proposition}
\noindent\textbf{Proof.} If $x\in\Sigma_{\mathcal{V}}g$, then $d_x\tilde{g}|_{V_{\alpha}}=0$, for a stratum $V_{\alpha}\in\mathcal{V}$ containing $x$ and an analytic extension $\tilde{g}$ of $g$ in a neighbourhood of $x.$ In $V_{\alpha}\subset\{f=0\}$, then $x=0,$ since $\Sigma_{\mathcal{V}}g\cap\{f=0\}=\{0\}.$ If $V_{\alpha}\subset X\setminus\{f=0\}$, $rk(d_x\tilde{f},d_x\tilde{g})\leqslant1,$ where $\tilde{f}$ is an analytic extension of $f$ in a neighbourhood of $x$. Hence, $x\in\Sigma(f,g)|_{V_{\alpha}}=\Sigma(f,g)|_{V_{\alpha}\setminus\{f=0\}}$, that is, $x\in\Gamma_{f,g}(\mathcal{V}).$ 

Furthermore, $\tilde{\Gamma}_{f,g}(\mathcal{V})$ is given by the components of $\Gamma_{f,g}(\mathcal{V})$ not contained in $\{g=0\},$ that is, $\tilde{\Gamma}_{f,g}(\mathcal{V})=\overline{\Gamma_{f,g}(\mathcal{V})\setminus\{g=0\}}\subset\Gamma_{f,g}(\mathcal{V}).$\qed

Using this proposition, we obtain the following useful information about the behaviour of $g$ with respect to $\mathcal{V}.$ 

\begin{corollary}\label{3.2}
If $\Gamma_{f,g}(\mathcal{V})$ is empty, then $g$ is prepolar at the origin with respect to the good stratification $\mathcal{V}$ of $X$ relative to $f.$
\end{corollary}
\noindent\textbf{Proof.} By Proposition \ref{3.1}, if $\Gamma_{f,g}(\mathcal{V})$ is empty, then $\Sigma_{\mathcal{V}}g$ is empty, that is, $g$ has no stratified critical point with respect to $\mathcal{V}.$\qed

Let $n_i$ be the number of stratified Morse critical points of a Morsefication of $g:X\cap f^{-1}(\delta)\cap B_{\epsilon}\rightarrow\mathbb{C}$ in $V_i\cap f^{-1}(\delta)\cap B_{\epsilon}$, for each $i\in\{1,\ldots,q\}.$ The next proposition uses the relative polar variety $\Gamma_{f,g}(\mathcal{V})$ for counting the number $n_i.$

\begin{proposition}\label{3.3}
If $\Gamma_{f,g}(\mathcal{V})$ is empty, then $n_i=0,$ for all $i\in\{1,\ldots,q\}.$
\end{proposition}
\noindent\textbf{Proof.} Let $V_i$ be a stratum of $\mathcal{V}$ and $x$ be a critical point of $g|_{V_i\cap f^{-1}(\delta)\cap B_{\epsilon}}.$ If $\tilde{f}$ and $\tilde{g}$ are analytic extensions of $f$ and $g$ in a neighbourhood of $x,$ respectively, then $x\in V_i\cap f^{-1}(\delta)\cap B_{\epsilon}$ and $rk(d_x\tilde{f},d_x\tilde{g})\leqslant1$, that is, $x\in(V_i\cap f^{-1}(\delta)\cap B_{\epsilon})\cap(\Sigma_{\mathcal{V}}f\cup\Sigma_{\mathcal{V}}g\cup\tilde{\Gamma}_{f,g}(\mathcal{V})).$ By Proposition 1.3 of \cite{Ms1}, $\Sigma_{\mathcal{V}}f\subset\{f=0\}.$ Therefore, by Proposition \ref{3.1},
\begin{center}
$\Sigma g|_{V_i\cap f^{-1}(\delta)\cap B_{\epsilon}}=V_i\cap f^{-1}(\delta)\cap B_{\epsilon}\cap(\Sigma g|_{V_i}\cup\tilde{\Gamma}_{f,g}(V_i))\subset V_i\cap f^{-1}(\delta)\cap B_{\epsilon}\cap\Gamma_{f,g}(V_i).$
\end{center}
Since $\Gamma_{f,g}(\mathcal{V})$ is empty, $\Sigma g|_{V_i\cap f^{-1}(\delta)\cap B_{\epsilon}}$ is empty which implies $n_i=0,$ for all $i\in\{1,\ldots,q\}.$\qed

In \cite{DG}, Dutertre and Grulha proved a Lê-Greuel type formula for the Brasselet number, making possible to count the number of stratified Morse critical points using Brasselet numbers. We apply their result to compute Brasselet numbers in the setting we already know the number of Morse critical points. 
First, let us show a more general result.

\begin{proposition}\label{3.4}
Let $\beta:X\rightarrow\mathbb{Z}$ be a constructible function with respect to the good stratification $\mathcal{V}$ of $X$ relative to $f$. If $\Gamma_{f,g}(\mathcal{V})$ is empty, then \begin{center}
$\chi(X\cap f^{-1}(\delta)\cap B_{\epsilon},\beta)=\chi(X\cap g^{-1}(0)\cap f^{-1}(\delta)\cap B_{\epsilon},\beta).$
\end{center}
\end{proposition} 
\noindent\textbf{Proof.} By Corollary \ref{3.2}, since $\Gamma_{f,g}(\mathcal{V})$ is empty, $g$ is prepolar at the origin with respect to $\mathcal{V}$ and, by Proposition 1.12 of \cite{Ms1}, $g$ is tractable at the origin with respect to $\mathcal{V}$. Then, by Theorem 4.2 of \cite{DG}, we obtain \begin{center}
$\chi(X\cap f^{-1}(\delta)\cap B_{\epsilon},\beta)=\chi(X\cap g^{-1}(0)\cap f^{-1}(\delta)\cap B_{\epsilon},\beta)+\sum_{i=1}^q(-1)^{d-1}n_i\eta(V_i,\beta),$
\end{center}
\noindent where $n_i$ is the number of stratified Morse critical points of a Morsefication of $g|_{V_i\cap f^{-1}(\delta)\cap B_{\epsilon}}$ appearing in $V_i\cap f^{-1}(\delta)\cap B_{\epsilon}.$ By Proposition \ref{3.3}, $n_i=0,$ for all $i\in\{1,\ldots,q\},$ and the equality holds.\qed

If the constructible function $\beta$ is the local Euler obstruction, we obtain a relation between Brasselet numbers.

\begin{corollary}\label{3.5}
If $X$ is equidimensional and $\Gamma_{f,g}(\mathcal{V})$ is empty, then $B_{f,X}(0)=B_{f,X^g}(0).$
\end{corollary}
\noindent\textbf{Proof.} By Theorem 4.4 of \cite{DG}, $B_{f,X}(0)=B_{f,X^g}(0)+(-1)^{d-1}n_q,$ where $n_q$ is the number of stratified Morse critical points of a Morsefication of $g|_{X_{reg}\cap f^{-1}(\delta)\cap B_{\epsilon}}$ appearing in $X_{reg}\cap f^{-1}(\delta)\cap B_{\epsilon}.$ Since $\Gamma_{f,g}(\mathcal{V})$ is empty, $n_q$ vanishes and the equality holds. \qed

When $f$ is a generic linear form on $\mathbb{C}^n, B_{f,X}(0)=Eu_X(0)$ and $B_{f,X^g}(0)=Eu_{X^g}(0).$ Therefore, the last corollary leads the following consequence.

\begin{corollary}
If $X$ is equidimensional and $\Gamma_{f,g}(\mathcal{V})$ is empty, then $Eu_X(0)=Eu_{X^g}(0).$
\end{corollary}

When both $f$ and $g$ have isolated singularity at the origin, Dutertre and Grulha proved several formulas about the Brasselet numbers of $f$ and $g$. Using these formulas and supposing that $\Gamma_{f,g}(\mathcal{V})$ is empty, we obtain further information about these numbers. If this is the case, then $g$ is prepolar at the origin with respect to the good stratification $\mathcal{V}$ of $X$ induced by $f$, given as a refinement of a Whitney stratification $\mathcal{W}=\{W_i\}$ of $X.$ Applying Corollary 6.1 of \cite{DG}, $f$ is prepolar a the origin with respect to the good stratification $\overline{\mathcal{V}}$ of $X$ induced by $f$, also given by a refinement of $\mathcal{W}.$ By Proposition 1.12 of \cite{Ms1}, we obtain that $\Gamma_{f,g}(\mathcal{V})=\tilde{\Gamma}_{f,g}(\mathcal{V})$ and $\Gamma_{g,f}(\overline{\mathcal{V}})=\tilde{\Gamma}_{g,f}(\overline{\mathcal{V}}).$ On the other hand, \begin{eqnarray*}
\tilde{\Gamma}_{f,g}(\mathcal{V})&=&\bigcup_{V_i\in\mathcal{V}}\overline{\Sigma(f,g)|_{V_i\setminus(\{f=0\}\cup\{g=0\})}}\\
&=&\bigcup_{W_i\in\mathcal{w}}\overline{\Sigma(f,g)|_{W_i\setminus(\{f=0\}\cup\{g=0\})}}\\
&=&\bigcup_{\overline{V_i}\in\overline{\mathcal{V}}}\overline{\Sigma(f,g)|_{\overline{V_i}\setminus(\{f=0\}\cup\{g=0\})}}\\
&=&\tilde{\Gamma}_{g,f}(\overline{\mathcal{V}}).
\end{eqnarray*} 

Hence, the four mentioned polar varieties are equal. We can then compute the following.

\begin{proposition}\label{3.7}
Let $\beta:X\rightarrow\mathbb{Z}$ be a constructible function with respect to the Whitney stratification $\mathcal{W}$ and $\mathcal{V}$ is the good stratification of $X$ induced by $f.$ If $\Gamma_{f,g}(\mathcal{V})$ is empty, then $\chi(X\cap f^{-1}(\delta)\cap B_{\epsilon},\beta)=\chi(X\cap g^{-1}(\delta)\cap B_{\epsilon},\beta).$
\end{proposition}
\noindent\textbf{Proof.} Since $\Gamma_{f,g}(\mathcal{V})$ is empty, by Corollary \ref{3.2}, $g$ is prepolar at the origin with respect to $\mathcal{V}.$ By Corollary 6.1 of \cite{DG}, $f$ is prepolar at the origin with respect to the good stratification $\overline{\mathcal{V}}$ induced by $g.$ By Theorem 6.4 of \cite{DG}, $$\chi(X\cap f^{-1}(\delta)\cap B_{\epsilon},\beta)=\chi(X\cap g^{-1}(\delta)\cap B_{\epsilon},\beta)+\sum_{i=1}^q(-1)^{d_i-1}(n_i-m_i)\eta(W_i,\beta),$$
\noindent where $d_i$ denotes the dimension of $W_i\in\mathcal{W}, n_i$ is the number of stratified Morse critical points of a Morsefication of $g|_{V_i\cap f^{-1}(\delta)\cap B_{\epsilon}}$ appearing in $V_i\cap f^{-1}(\delta)\cap B_{\epsilon}$ and $m_i$ is the number of stratified Morse critical points of a Morsefication of $f|_{\overline{V_i}\cap g^{-1}(\delta)\cap B_{\epsilon}}$ appearing in $\overline{V_i}\cap g^{-1}(\delta)\cap B_{\epsilon}.$ By Proposition \ref{3.3}, since $\Gamma_{f,g}(\mathcal{V})$ is empty, $m_i=n_i=0,$ for all $i\in\{1,\ldots,q\},$ and the equality is proved. \qed

In the case that $\beta$ is the local Euler obstruction, we obtain a relation between Brasselet numbers. 

\begin{corollary}
If $X$ is equidimensional and $\Gamma_{f,g}(\mathcal{V})$ is empty, $B_{f,X}(0)=B_{g,X}(0).$
\end{corollary}

\noindent\textbf{Proof.} Since $\Gamma_{f,g}(\mathcal{V})$ is empty, $g$ (resp. $f$) is prepolar at the origin with respect to the good stratification of $f$ (resp. $g$). By Corollary 6.5 of \cite{DG}, $$B_{f,X}(0)=B_{g,X}(0)+(-1)^{d-1}(n_q-m_q),$$ where $n_q$ (resp. $m_q$) is the number of Morse critical points of a Morsefication of $g|_{X_{reg}\cap f^{-1}(\delta)\cap B_{\epsilon}}$ (resp. $f|_{X_{reg}\cap g^{-1}(\delta)\cap B_{\epsilon}}$)  appearing in $X_{reg}\cap f^{-1}(\delta)\cap B_{\epsilon}$ (resp. $X_{reg}\cap g^{-1}(\delta)\cap B_{\epsilon}$). Using again that $\Gamma_{f,g}(\mathcal{V})$ is empty, we have that $n_q=m_q=0$ and the equality holds.\qed

 \end{document}